\documentclass[12pt]{article}

\usepackage{amsmath}[1996/11/01]
\usepackage{amssymb,amsthm,amsfonts,latexsym}
\usepackage{amscd, makeidx}

\def\beql#1#2\eeql{\begin{equation}\label{#1}#2\end{equation}}

\newtheorem{thm}{Theorem}
\newtheorem{theorem}[thm]{Theorem}
\newtheorem{kor}[thm]{Corollary}
\newtheorem{cor}[thm]{Corollary}

\newtheorem{lemma}[thm]{Lemma}

\theoremstyle{definition}

\newcommand{\bew}{\noindent\underline{Proof.}\ }
\newcommand{\eb}{\phantom{zzz}\hfill{$\square $}\smallskip}

\DeclareMathOperator{\GL}{GL}

\DeclareMathOperator{\Aut}{Aut}
\DeclareMathOperator{\Tr}{Tr}
\DeclareMathOperator{\tr}{trace}
\DeclareMathOperator{\trace}{trace}

\newcommand{\C}{\mathbb C}
\newcommand{\R}{\mathbb R}
\newcommand{\Q}{\mathbb Q}

\newcommand{\Z}{\mathbb Z}

\newcommand{\N}{\mathbb N}

\begin{document}

\LARGE
\begin{center}
{\bf Hermitian modular forms congruent to 1 modulo $p$.}
\end{center}

\normalsize
\begin{center}
Michael Hentschel 
\\Lehrstuhl A f\"ur Mathematik, RWTH Aachen University, 52056 Aachen, Germany, hentschel@mathA.rwth-aachen.de \\
Gabriele Nebe
\\Lehrstuhl D f\"ur Mathematik, RWTH Aachen University, 52056 Aachen, Germany, nebe@math.rwth-aachen.de
\end{center}

\begin{abstract}
For any natural number $\ell $ and any prime $p\equiv 1 \pmod{4}$ 
not dividing $\ell $ 
there is a Hermitian modular form of arbitrary genus $n$
over $L:=\Q [\sqrt{-\ell }]$ that
is congruent to $1$ modulo $p$ which is a Hermitian theta series of an
$O_L$-lattice of rank $p-1$ 
admitting a fixed point free automorphism of order $p$.
It is shown that also for non-free lattices such theta series are
modular forms.
\end{abstract}

\section{Introduction.}

The purpose of the present note is to generalize the construction of 
Siegel modular forms that are congruent to 1 modulo a suitable prime $p$ 
given in \cite{Boe} to the case of Hermitian modular forms over 
$L:=\Q [\sqrt{-\ell }]$. For $\ell = 1$ and $\ell = 3$ this was done in
\cite{Nag}, in fact we use the same strategy by constructing an
even unimodular lattice $\Lambda $ as an ideal lattice in 
$K:= L [\zeta _p]$ for any prime $p\equiv 1 \pmod{4}$ not dividing $\ell $.
The existence of $\Lambda $ essentially follows from class field theory 
and is predicted by \cite[Th\'eor\`eme 2.3, Proposition 3.1 (1)]{BaMa}
(see also \cite[Corollary 2]{Ba}). 
Since the ring of integers $O_L$ is in general not a principal ideal 
domain the lattice $\Lambda $ is not necessarily a free $O_L$-module. 
We are not aware of an explicit statement in the literature that 
the genus $n$ Hermitian theta series $\theta ^{(n)}(\Lambda )$
of such a lattice $\Lambda $ is 
a modular form for the full modular group. Therefore 
the first section sketches
a proof. In fact the proofs in the literature never seriously use the fact that 
the lattice is a free $O_L$-module. 
The next section applies the results of \cite{BaMa} and \cite{Ba}
to the special case of the 
field $K=\Q [\sqrt{-\ell },\zeta _p ]$ and proves the existence of 
a Hermitian $O_K$-lattice $\Lambda _h$ that is an even unimodular $\Z $-lattice 
(with respect to the trace of the Hermitian form). 
The invariance under $O_K$ yields both, a Hermitian $O_L$-module structure 
on $\Lambda _h$ and an $O_L$-linear automorphism (the multiplication by 
the primitive $p$-th root of unity $\zeta _p \in O_K$) of order
$p$ acting fixed point freely on $\Lambda _h \setminus \{0 \}$. 
Therefore all but the first coefficient in $\theta ^{(n)} (\Lambda _h)$ 
are multiples of $p$ yielding the desired Hermitian modular form.

\section{Hermitian theta-series are Hermitian modular forms.}

Let $\ell \in \N $ such that $-\ell $ is a fundamental discriminant
(which means that 
either $\ell \equiv -1 \pmod{4} $  is square-free 
or
$\ell =4m$, where $m\equiv 2$ or $1 \pmod{4}$  is square-free).
Let $L:= \Q [\sqrt{-\ell }]$ be the imaginary quadratic number 
field of discriminant $-\ell $, with ring of integers $O_L$ 
and  inverse different 
$$O_L^* := \{ a \in L \mid \Tr _{L/\Q}  (a O_L) \subset \Z \} = \sqrt{-\ell }^{-1} O_L .$$

 Let $(V,h)$ be a finite dimensional 
positive definite Hermitian vector space over $L$.
This section extends the results in \cite{CohRes} to not necessarily free 
Hermitian $O_L$-lattices in $(V,h)$. 
Note that we use a different scaling for the Hermitian form 
resulting in the additional factor of 2 in the definition of the
Hermitian Siegel theta series below.
It is already stated in \cite{CohRes} that the authors restrict to free lattices
``for convenience'' and that the same results hold in the more general context.
The {\bf full modular group} 
$$\Gamma _n := \langle 
\left( \begin{array}{cc} I_n & B \\ 0 & I_n \end{array} \right) ,\ 
\left( \begin{array}{cc} U & 0 \\ 0 & U^{-1} \end{array} \right) ,\ 
\left( \begin{array}{cc} 0 & -I_n \\  I_n & 0 \end{array} \right)  \mid
B \in O_L^{n\times n} \mbox{ Hermitian }, U \in \GL_n(O_L) \rangle $$
(see \cite{Klingen}, \cite[Anhang V]{Freitag}, 
\cite{SkriptKrieg})  for the proof that these matrices really generate) 
acts on the Hermitian half space 
by 
$$ Z \mapsto Z+B,\ Z \mapsto \overline{U}^{t} Z U ,\ Z \mapsto -Z^{-1} $$
for the respective generators.

\begin{theorem}\label{nonfree}
Let $\Lambda _h \leq (V,h)$ be an $O_L$-lattice such that 
the $O_L$-dual lattice
$$\Lambda _h ^{*} := \{ v\in V \mid h(v,\Lambda _h ) \subset O_L \}  
= \sqrt{-\ell } \Lambda _h  =  (O_L^*)^{-1} \Lambda _h.$$ 
Then its Hermitian theta series
$$\theta ^{(n)} (\Lambda _h)(Z) := \sum _{(x_1,\ldots , x_n)\in \Lambda _h ^n}
\exp (2\pi i \tr (h(x_i,x_j) Z))$$
is a Hermitian modular form for the full modular group $\Gamma _n$.
Here $\tr : L^{n\times n} \to \Q $ denotes the composition of the
matrix trace with the trace of $L$ over $\Q $. 
\end{theorem}

\bew
For $x:=(x_1,\ldots ,x_n) \in \Lambda _h^n$ the Hermitian matrix 
$H:=H_x:= (h(x_i,x_j)) \in  (O_L^*)^{n\times n}$, so 
for any Hermitian matrix $B \in O_L^{n\times n }$ the trace 
$\tr (HB) $ is in $\Z $. 
This shows the invariance of $\theta ^{(n)} (\Lambda _h)$ under $Z\mapsto Z+B$.
Similarly 
$$\tr (H _x \overline{U}^t Z U ) = \tr (U H _x \overline{U^t} Z) = 
\tr (H _{xU} Z) $$ 
so the transformation $Z\mapsto \overline{U}^t Z U $ for $U\in \GL_n(O_L)$ 
just changes the order of summation in $\theta ^{(n)} (\Lambda _h)$.
It remains to prove the theta-transformation formula 
$$ (\star ) \ \  \theta ^{(n)}(\Lambda _h) (-Z^{-1}) = \det(Z/i)^d  \theta ^{(n)}(\Lambda _h) (Z) $$
also 
for non-free $O_L$-lattices $\Lambda _h$ of dimension $d$ 
that satisfy $\Lambda _h^* = (O_L^*)^{-1} \Lambda _h$.
But Poisson summation only depends on the abelian group structure, 
not on the underlying module, so the usual proof 
(see for instance \cite[p. 111]{Krieg})
can be adopted to the situation here 
(for details we refer to \cite{thesis}):
Using the Identity Theorem, it suffices to prove $(\star )$ for $Z = iY $, 
$Y$ Hermitian positive definite. 
Let $\varphi : \R^{2dn} \to \C^{d \times n } $ be the obvious isomorphism
and consider $\Lambda _h^n$ as a lattice $\tilde{\Lambda }$ 
in $\C^{d \times n } $ choosing
coordinates with respect to an orthonormal basis of $(\C^d,h)$.
Then there is some $F \in \R ^{2dn\times 2dn} $ such that 
$\tilde{\Lambda } = \varphi (F \Z ^{2dn} ) $ and  $H_{\varphi (x)} = \overline{\varphi (Fx) }^{tr} \varphi(Fx)$. Then
$$ \theta ^{(n)}(\Lambda _h) (iY) = \sum _{g\in \Z ^{2dn}} \psi (g) $$
where 
$$\psi : \R ^{2dn} \to \C , 
x \mapsto \exp (-2 \pi \trace (\overline{\varphi (Fx) }^{tr} \varphi(Fx) Y)) .$$
The condition $\Lambda _h^* = (O_L^*)^{-1} \Lambda _h$ implies that 
$|\det (F) | = 1$ and we can apply the usual Poisson summation to get the result 
as in \cite[pp. 110-112]{Krieg}.
\eb

\section{Congruences of Hermitian theta-series.}

Let $p$ be a prime $p\equiv 1 \pmod{4}$ such that $\ell $ is not a multiple of $p$. 
This section constructs a Hermitian $O_L$-lattice
$(\Lambda ,h) $ of rank $p-1$ admitting an automorphism of order $p$
such that the $\Z $-lattice 
$(\Lambda, \Tr _{L/\Q}  (h) )$ is a positive definite even unimodular lattice. 
The existence of such a lattice follows from the much more general 
result \cite[Th\'eor\`eme 2.3]{BaMa} together with 
\cite[Proposition 3.1]{BaMa} which are
based on Artin's reciprocity law in global class field theory
(see \cite[Theorem (V.3.5)]{Milne}).
For our special case it is however more convenient to use 
\cite[Corollary 2]{Ba}, which is essentially a consequence of \cite[Th\'eor\`eme 2.3]{BaMa}. 

To this aim we consider the number field 
$K = \Q[\sqrt{-\ell }][\zeta _p ]=LM$ with $M=\Q[\zeta _p]$,
 where $\zeta _p = \exp (\frac{2\pi i}{p} )$
is a primitive $p$-th root of unity. 
Then $K$ is an abelian number field of degree $2(p-1)$ over $\Q $ which is 
a multiple of 8.
The field $K$ is totally complex and admits an involution 
$\overline{\phantom{x}}$, the complex conjugation, with 
fixed field $F$ the totally real subfield of $K$.

The following lemma is well known.

\begin{lemma}\label{unram}
$K/F$ is unramified at all finite primes.
\end{lemma}

\bew
The discriminant $d_{K/F}$ of $K/F$ divides the discriminant of 
any $F$-basis of $K$ that consists of integral elements. 
For $B_1= (1,\sqrt{-\ell} )$ 
one finds $d_{B_1} = \det (\Tr _{K/F} (b_i b_j ))  = 
-4\ell $ and for 
$B_2 = (1,\zeta _p)$ one get $d_{B_2} = \zeta _p^{-2} ( \zeta _p^2-1)^2 $ 
which generates an ideal of norm $p^2$ in $F$.
Since $p$ is an odd prime not dividing $\ell $, the $\gcd $ of these two 
discriminants is $1$ and hence $d_{K/F} = 1$ which implies the lemma.
\eb

Since all real embeddings of $F$ extend to complex embeddings of $K$ and 
$[K:\Q] = 2 (p-1) \equiv 0 \pmod{8} $ \cite[Corollary 2]{Ba} yields the existence 
of a fractional $O_K$-ideal ${\cal A}$ in $K$ and a totally positive element 
$d \in F$ such that the $O_K$-module ${\cal A}$ together with the 
symmetric integral bilinear form 
$$ b_d : {\cal A} \times {\cal A} \to \Z , \ (x,y) \mapsto \trace_{K/\Q} (dx\overline{y}) $$
is an even unimodular $\Z$-lattice $\Lambda := ({\cal A}, b_d) $.
This means that $b_d(x,x) \in 2\Z $ for all $x \in {\cal A} $ and 
$$\Lambda ^{\#} :=
 \{ x\in K \mid b_d(x,y) \in \Z \mbox{ for all } y\in {\cal A} \}  = \Lambda .$$

\begin{cor}\label{uni}
The $O_L$-lattice $\Lambda _h:=({\cal A},h(x,y) := \Tr _{K/L} (d x \overline{y })) $
is a Hermitian $O_L$-lattice with automorphism 
$x\mapsto \zeta _p x$ of order $p$ such that $\Lambda _h^* = (O_L^*)^{-1} \Lambda _h$. 
\end{cor}

\bew
Since ${\cal A}$ is an ideal of $K$, the 
 multiplication by $\zeta _p \in O_K$ preserves the lattice ${\cal A}$.
It also respects  
the Hermitian form $h$, because 
$$h(\zeta_px,\zeta_py) = \Tr _{K/L} (d \zeta_p x \overline{\zeta_p y }))
= \Tr _{K/L} (d \zeta_p \zeta _p^{-1} x \overline{ y }))= h(x,y) .$$
The fact that  $\Lambda _h^* = (O_L^*)^{-1} \Lambda _h$ follows from the unimodularity of the
integral lattice $\Lambda $:
For $y\in K$ we obtain 
$$b_d(x,y ) = \trace _{L/\Q}(h(x,y)) \in \Z \mbox{ for all } x \in {\cal A} 
\Leftrightarrow h(x,y) \in O_L^*   \mbox{ for all } x \in {\cal A} $$
using the fact that ${\cal A}$ is an $O_L$-module and $h$ is Hermitian over $O_L$.
Hence   $\Lambda _h^* =  (O_L^*)^{-1} \Lambda ^{\#} = (O_L^*)^{-1} \Lambda _h$.
\eb

Together this implies the existence of a Hermitian modular 
form of weight $p-1$ that is congruent to 1 modulo p for 
more general imaginary quadratic number fields than those treated in
\cite{Nag}:

\begin{theorem}
Let $L = \Q[\sqrt{-\ell }]$ be an imaginary quadratic number field
($-\ell $ a fundamental discriminant) 
and let $p$ be a prime $p\equiv 1 \pmod{4}$ not dividing $\ell $.
Then for arbitrary genus $n \geq 1$ there is a Hermitian modular form 
$$F_{p-1}^{(n)} \in M_{p-1}(SU_n(O_L)) $$ 
for the full modular group over the ring of integers $O_L$ of $L$ such that 
$$F_{p-1}^{(n)} \equiv 1 \pmod{p} .$$
\end{theorem}

\bew
Corollary \ref{uni}  constructs 
a Hermitian $O_L$-lattice $\Lambda _h$ of rank 
$p-1$ admitting an automorphism of order $p$ (which necessarily acts fixed point
freely) such that $\Lambda _h^* = (O_L^*)^{-1} \Lambda _h$. 
By  Theorem \ref{nonfree} its Siegel theta series 
is a Hermitian modular form for the full modular group.
Since $\Lambda _h$ admits a fixed point free automorphism of order $p$,
all the representation numbers 
$$R_A := | \{ (x_1,\ldots , x_n)\in \Lambda ^n \mid (h(x_i,x_j)) = A \} |$$
for any non-zero Hermitian matrix $A \in L^{n\times n }$ are 
multiples of $p$ and hence 
$$ F_{p-1}^{(n)} := \theta ^{(n)}_{\Lambda _h}  \equiv 1 \pmod{p} $$ 
provides the desired Hermitian modular form.
\eb

Since the root lattice $E_8$ is the unique even unimodular lattice of 
dimension 8, we obtain the following corollary.

\begin{kor}
Let $\ell \in \N $ be not a multiple of $5$.
Then the root lattice $E_8$ has a Hermitian structure as a 
lattice $\Lambda _h$ over 
the ring of integers of $\Q [\sqrt{-\ell }]$ such that 
$\Aut (\Lambda _h)$ contains an element of order $5$. 
\end{kor}



\begin{thebibliography}{1}
\bibitem{Ba}
E. Bayer-Fluckiger, 
Determinants of integral ideal lattices and automorphisms of given
characteristic polynomial. J. Algebra, 257 (2002), 215-221.
\bibitem{BaMa}
E. Bayer-Fluckiger, J. Martinet, 
Formes quadratiques lie\'es aux alg\`ebres semi-simples.
{\em J. reine angew. Math.} {\bf 451} (1994) 51-69
\bibitem{Boe} S. Boecherer, S. Nagaoka, On mod p properties of 
Siegel modular forms. {\em Math. Ann.} {\bf 338} (2007) 421-433
\bibitem{CohRes}  D.M. Cohen, H.L. Resnikoff, 
Hermitian quadratic forms and hermitian modular forms. 
{\em Pacific J. Math.} {\bf 76} 329-337 (1978)
\bibitem{Freitag} E. Freitag, {\em Siegelsche Modulfunktionen.} 
Grundlehren der mathematischen Wissenschaften {\bf 254} Springer (1983)
\bibitem{thesis} M. Hentschel, {\em On Hermitian theta 
series and modular forms.} PhD thesis, RWTH Aachen 2009. 
\bibitem{Klingen}  H. Klingen, 
Bemerkung \"uber Kongruenzuntergruppen der Modulgruppe n-ten Grades.
 {\em Arch. Math.} {\bf 10} 113-122 (1959)
\bibitem{SkriptKrieg} A. Krieg, 
Siegelsche Modulformen. Skript RWTH Aachen (2007)
\bibitem{Krieg} A. Krieg, 
{\em Modular forms on half-spaces of quaternions.} Lecture Notes in Mathematics, 1143. Springer-Verlag, Berlin, 1985.
\bibitem{Milne} J.S. Milne, {\em Class field theory.} Lecture Notes available 
via http://www.jmilne.org/math/
\bibitem{Nag}
T. Kikuka, S. Nagaota, Congruence properties of Hermitian modular forms.
Preprint 2008. 
\end{thebibliography}
\end{document}